\documentclass[10pt]{article}

\makeatletter
\@addtoreset{equation}{section}
\makeatother

\newtheorem{remark}{Remark}

\newtheorem{example}{Example}

\usepackage{epsfig}
\usepackage{amssymb}
\usepackage{amsmath}
\usepackage{graphicx}

\def\tfrac#1#2{{{\lower.6ex
\hbox{$\scriptstyle#1$}}\over
{\raise.7ex
\hbox{$\scriptstyle#2$}}}}

\def\Frac#1#2{\frac{\displaystyle{#1}}{\displaystyle{#2}}}

\def\dsp#1{\displaystyle#1}

\def\protectbold#1{\protect{\boldmath{$#1$}}}

\def\calC{{\cal C}}

\def\sign{{\rm sign}}

\def\erfc{{\rm erfc}}

\begin{document}

\title{The asymptotic and numerical inversion of the Marcum $Q-$function}

\author{Amparo Gil\\
Departamento de Matem\'atica Aplicada y Ciencias \\
de la Computaci\'on ETSI Caminos. \\
Universidad de Cantabria. 39005-Santander, Spain.\\
\and
Javier Segura\\
        Departamento de Matem\'aticas, Estad\'{\i}stica y 
        Computaci\'on,\\
        Universidad de Cantabria, 39005 Santander, Spain.\\   
\and
Nico M. Temme\\
IAA, 1391 VD 18, Abcoude, The Netherlands\footnote{Former address: CWI, Science Park 123, 1098 XG Amsterdam, The Netherlands}\\ \\
{ \small
  e-mail: {\tt
    amparo.gil@unican.es,
    javier.segura@unican.es, 
    nico.temme@cwi.nl}}
    }

\maketitle
\begin{abstract}
 The generalized Marcum  functions appear in
problems of technical and scientific areas such as, for example, radar detection and communications. In mathematical
statistics and probability theory these functions are called the noncentral gamma or the 
noncentral chi-squared cumulative distribution functions.
In this paper we describe a new asymptotic method for inverting the generalized Marcum $Q-$function and 
for the complementary Marcum $P-$function. Also, we show how monotonicity and convexity properties of these functions
can be used to find initial values for reliable Newton or secant methods to invert the function. We present details of 
numerical computations that show the reliability of the asymptotic approximations.

\end{abstract}

\vskip 0.8cm \noindent
{\small
2000 Mathematics Subject Classification:
33C10, 33B20, 41A60, 65D20 .
\par\noindent
Keywords \& Phrases:
Marcum $Q-$function, noncentral gamma distribution, noncentral $\chi^2$-distribution, 
incomplete gamma functions,  asymptotic expansions, numerical inversion.
}

\section{Introduction}\label{intro}
We define the generalized Marcum $Q-$function by using the integral representation
\begin{equation}
\label{eq:intro01}
Q_{\mu} (x,y)=\displaystyle x^{\frac12 (1-\mu)} \int_y^{+\infty} t^{\frac12 (\mu -1)} e^{-t-x} I_{\mu -1} \left(2\sqrt{xt}\right) \,dt,
\end{equation}
where $\mu >0$ and $I_\mu(z)$ is the modified Bessel function. 
We also use the complementary function
\begin{equation}
\label{eq:intro02}
P_{\mu} (x,y)=\displaystyle x^{\frac12 (1-\mu)}\int_0^{y} t^{\frac12 (\mu -1)} e^{-t-x} I_{\mu -1} \left(2\sqrt{xt}\right) \,dt,
\end{equation}
and the complementary relation reads
\begin{equation}\label{eq:intro03}
P_{\mu}(x,y)+Q_{\mu} (x,y)=1.
\end{equation}

There are other notations for the generalized Marcum function in the literature. Among them, probably the most
 popular is the following
\begin{equation}\label{eq:intro04}
\widetilde{Q}_{\mu} (\alpha,\beta)=\alpha^{1-\mu}\int_\beta^{+\infty}t^{\mu} e^{-(t^2+\alpha^2)/2} I_{\mu -1}(\alpha t) dt
\end{equation}
where we have added a tilde in the definition to distinguish it from the definition we are using (\ref{eq:intro01}). For
$\mu=1$ the definitions coincide with the original definition of the Marcum $Q-$function \cite{Marcum:1960:AST}. The relation with the notation we use is simple:
\begin{equation}\label{eq:intro05}
Q_{\mu}(x,y)=\widetilde{Q}_{\mu} (\sqrt{2x},\sqrt{2y}),
\end{equation}
and similarly for the $P$ function.

The generalized Marcum $Q-$function is an important function used in many applications in science and engineering, 
and notably in radar detection and communications. 
These functions also occur in statistics and probability theory, where they are called noncentral chi-squared or noncentral 
gamma cumulative distributions. 
Noncentral distributions play an important role in statistics because they arise in the power analysis of
statistical tests.
The central gamma cumulative distribution is in fact the incomplete gamma function, and the relation of the Marcum functions with
 the incomplete gamma functions is given in \S\ref{sec:props}. 
 For references on these application areas we refer to our recent publication \cite{Gil:2013:CMQ}, in which we have described 
reliable numerical algorithms for computing $P_{\mu}(x,y)$ and $Q_{\mu} (x,y)$ for a
 wide range of the positive real parameters $\mu, x, y$ ($\mu\ge 1$).

In the present paper we describe a new method for inverting the generalized Marcum functions for large values of the parameter $\mu$. For this we use an asymptotic representation, and essential steps in the inversion algorithm are based on the inversion of the complementary error function. We have used the same approach in our recent paper \cite{Gil:2012:IGR} for numerical and asymptotic inversion  algorithms  for the incomplete gamma ratios. Also, we show how monotonicity and convexity properties of these functions
can be used to find initial values for reliable Newton or secant methods to invert the functions. A combination of both asymptotic
and numerical methods gives an efficient method of inversion for positive real values of the variables (with $\mu\ge 1$).

\section{Properties of Marcum functions}\label{sec:props}

Considering the Maclaurin series for the modified Bessel function and integrating term by term in (\ref{eq:intro01})
we obtain the series expansions
\begin{equation}\label{eq:def02}
\begin{array}{l}
\dsp{P_{\mu} (x,y)=e^{-x}\sum_{n=0}^{\infty} \frac{x^n}{n!}  P_{\mu +n}(y),}\\[8pt]
\dsp{Q_{\mu} (x,y)=e^{-x}\sum_{n=0}^{\infty} \frac{x^n}{n!}  Q_{\mu +n}(y).}
\end{array}
\end{equation}
These expansions are in terms of the incomplete gamma function ratios defined by
\begin{equation}\label{eq:def03}
P_{\mu} (x)=\frac{\gamma (\mu,x)}{\Gamma(\mu)},\quad
Q_{\mu} (x)=\frac{\Gamma (\mu,x)}{\Gamma(\mu)},
\end{equation}
where for $\Re\mu>0$ the standard incomplete gamma functions are defined by
\begin{equation}\label{eq:def04}
\gamma (\mu,x)=\int_{0}^x t^{\mu-1} e^{-t} \,dt,\quad
\Gamma (\mu,x)=\int_{x}^{+\infty} t^{\mu-1} e^{-t} \,dt.
\end{equation}
We have the complementary relation $P_\mu(x)+Q_\mu(x)=1$. 
Note that for the incomplete gamma function ratios $P_{\mu} (y)$ and $Q_{\mu} (y)$
appearing in (\ref{eq:def02}) algorithms are given in \cite{Gil:2012:IGR}.

The series expansion for $P_{\mu}(x,y)$ and $Q_{\mu}(x,y)$ given in \eqref{eq:def02} provide the standard
definition for the noncentral gamma distribution in statistics: if $Y$ is a a random
variable with a noncentral gamma distribution with parameter $a>0$ and noncentral parameter
$\lambda \ge 0$, then the  lower and upper tail probabilities
for the distribution function of $Y$ are given by \cite{pat:1949:NCS}

\begin{equation}\label{eq:defst}
\begin{array}{l}
\dsp{F(y|a,\lambda)=\mbox{Prob}(  Y <y) =\sum_{n=0}^{\infty}p(n|\lambda) F(y|a+n),}\\[8pt]
\dsp{F^c(y|a,\lambda)=\mbox{Prob}(  Y >y) =\sum_{n=0}^{\infty}p(n|\lambda) F^c(y|a+n),}
\end{array}
\end{equation}
respectively. Here,  $p(n|\lambda) =\Frac{e^{-\lambda} \lambda^n}{n!}$ is the point probability of a Poisson distribution, and $F(y|a)$, $F^c(y|a)$ are the lower and upper tail
probabilities, respectively, of a central gamma distribution function with parameter $a$; i.e., $F(y|a)$, $F^c(y|a)$ are the incomplete gamma function ratios $P_{a}(y)$ and $Q_a(y)$, respectively. In statistical terminology, it is said that the noncentral gamma distribution is a mixture of central gamma distributions with Poisson weights. 
From the noncentral gamma, the noncentral chi-squared distribution function is derived very easily: if 
the random variable $Y$ has a noncentral gamma distribution with parameters $a$ and non-centrality parameter 
$\lambda$, then $X=2Y$ has a noncentral chi-squared distribution with parameter $n=2a$ and with non centrality
parameter $2\lambda$.    

Particular values can be obtained from (\ref{eq:intro01}) and (\ref{eq:def02}):
\begin{equation}\label{eq:def07}
\begin{array}{ll}
Q_{\mu} (x, 0)=1,&Q_{\mu}(x ,+\infty)=0,\\[8pt]
Q_{\mu} (0,y)=Q_\mu(y),\quad& Q_{\mu}(+\infty,y)=1,\\[8pt]
Q_{+\infty} (x ,y )=1,&
\end{array}
\end{equation}
and similar complementary relations for $P_\mu(x,y)$.

As will follow from the relations given later (and in less detail from the relations in \eqref{eq:def07}), the transition in the $(x,y)$ quarter plane from small values of $Q_{\mu} (x,y)$ to values close to unity occurs for large values of $\mu, x,y $ across the  line $y=x+\mu$. Above this line in the $(x,y)$ quarter plane $Q_{\mu} (x,y)$ is smaller than $P_{\mu} (x,y)$ and below this   line the complementary function $P_{\mu} (x,y)$ is the smaller one. For more details on this transition line we refer to Example~\ref{ex:trans} in \S\ref{sec:asinvy}.

As we see next, this transition line is close to the inflection points of the graphs of Marcum functions both as a function of
$x$ and $y$.

\subsection{Monotonicity and convexity properties}

For the inversion process it is important to describe the monotonicity and convexity with respect to $x$ or $y$, with $\mu$ fixed.

First we mention the recurrence relations
\begin{equation}\label{eq:def05}
\begin{array}{ll}
\dsp{Q_{\mu+1}(x,y)=Q_{\mu} (x,y)+\left(\Frac{y}{x}\right)^{\mu/2} e^{-x-y} I_{\mu} (2\sqrt{xy})},\\[8pt]
\dsp{P_{\mu+1}(x,y)=P_{\mu} (x,y)-\left(\Frac{y}{x}\right)^{\mu/2} e^{-x-y} I_{\mu} (2\sqrt{xy})},
\end{array}\end{equation}
and the related three-term recurrence
\begin{equation}
\label{TTRR}
\begin{array}{l}
Q_{\mu+1}(x,y)-\left(1+c_{\mu}(x,y)\right)Q_{\mu}(x,y)+c_{\mu}Q_{\mu-1}(x,y)=0,\\[8pt]
c_{\mu}(x,y)=\sqrt{y/x}I_{\mu}(2\sqrt{xy})/I_{\mu -1}(2\sqrt{xy}),
\end{array}
\end{equation}
which is also satisfied by $P_{\mu}(x,y)$.

Taking the derivative with respect to $y$ in (\ref{eq:intro01}) and using (\ref{eq:def05}) we have
\begin{equation}\label{eq:def08}
\Frac{\partial Q_\mu(x,y)}{\partial y}=Q_{\mu-1}(x,y)-Q_{\mu}(x,y),
\end{equation}
and similarly
\begin{equation}\label{eq:def09}
\Frac{\partial Q_\mu(x,y)}{\partial x}=Q_{\mu+1}(x,y)-Q_{\mu}(x,y) .
\end{equation}
By using the relations in \eqref{eq:def05} it follows that
\begin{equation}\label{eq:def10}
\Frac{\partial Q_\mu(x,y)}{\partial x}=-\Frac{\partial Q_{\mu+1}(x,y)}{\partial y}= 
\left(\Frac{y}{x}\right)^{\mu/2} e^{-x-y} I_{\mu} \left(2\sqrt{xy}\right),
\end{equation}
and we see that $Q_\mu (x,y)$ ($P_\mu (x,y)$) is an increasing (decreasing) 
function of $x$ and a decreasing (increasing) function of $y$. With respect to $\mu$, $Q_\mu(x,y)$ is increasing and $P_\mu(x,y)$ is decreasing.

Regarding the convexity properties, using (\ref{TTRR}), (\ref{eq:def08}) and  (\ref{eq:def10}) we obtain
\begin{equation}
\label{dedps}
\begin{array}{ll}
\Frac{\partial^2 Q_{\mu}(x,y)}{\partial x^2}=(c_{\mu+1}(x,y)-1)(Q_{\mu+1}(x,y)-Q_{\mu}(x,y))
,\\[8pt]
\Frac{\partial^2 Q_{\mu}(x,y)}{\partial y^2}=\Frac{\partial^2 Q_{\mu-2}(x,y)}{\partial x^2},
\end{array}
\end{equation}
with $c_\mu(x,y)$ as defined in (\ref{TTRR}).

It is easy to prove, using similar ideas as in \cite{Segura:2011:BRM}, that $c_{\nu}(x,y)$ is decreasing as a 
function of $x$ and increasing as a function of $y$. On the other hand, we have $c_{\nu}(0,y)=y/\nu$. 
Therefore, if $\mu +1>y$ (both fixed) 
 we have $c_{\mu +1}(x,y)<1$ for all $x>0$, because $c_{\mu +1}(x,y)$ decreases as a function of $x$. Therefore
$$
\Frac{\partial^2 Q_{\mu}(x,y)}{\partial x^2}<0 \quad \mbox{ if } \quad  y<\mu +1,\ x>0.
$$

For $y>\mu+1$ the convexity may change as a function of $x$ (only once, because $c_{\mu+1}$ is monotonic). 
The inflection points with respect to $x$ and $y$ can be estimated quite sharply, as shown in 
\cite{Seg:2013:CPB} (by estimating the values for which $c_{\mu+1}(x,y)=1$). 
It can be proved that if $\mu > 0$ then
\begin{equation}
\begin{array}{l}
\Frac{\partial^2 Q_{\mu}(x,y)}{\partial x^2}<0 \quad \mbox{ if } \quad x>y-\mu-\dsp{\tfrac12},\\[8pt]
\Frac{\partial^2 Q_{\mu}(x,y)}{\partial x^2}>0 \quad \mbox{ if } \quad x<y-\mu-1,
\end{array}
\end{equation}
and if $\mu\ge \frac32$ then
\begin{equation}
\label{conve_y}
\begin{array}{l}
\Frac{\partial^2 Q_{\mu}(x,y)}{\partial y^2}>0 \quad \mbox{ if } \quad y>x+\mu-1,\\[8pt]
\Frac{\partial^2 Q_{\mu}(x,y)}{\partial y^2}<0 \quad \mbox{ if } \quad y<x+\mu-\dsp{\tfrac32}.
\end{array}
\end{equation}

As we discuss later, the monotonicity and convexity properties are useful for finding initial values for
the Newton or the secant methods, and they provide a reliable numerical inversion method for the Marcum functions.

\section{Inversion of Marcum's $Q-$function}\label{sec:MQinv}

In statistics, the inversion of the Marcum functions
refers to two different kinds of problems:

\begin{description}
\item{a)} Computation of the quantiles of the distribution; i.e., compute the $p-$quantile 
$y$ for which $F(y|a,\lambda)=p$.

\item{b)} Computation of the non-centrality parameter of the distribution given the lower or upper
tail probability; i.e., for given $y,\,a,$ and $p$ obtain $\lambda$ such that
$F(y|a,\lambda)=p$ or $F^c(y|a,\lambda)=p$.  
\end{description}

Newton algorithms for solving both problems have been proposed in the literature; see for example, \cite{ding:1999:NEW}
for problem a) and \cite{knu:1996:GAM}. 
No special study when the parameter $a$ is large is discussed in any of these references.

To discuss the inversion of the Marcum $Q-$function, we follow the process described by Helstrom \cite{Helstrom:1998:AIM}, where
the inversions are  linked to a  specific problem in radiometry.
In this reference, the inversion of the $Q-$function is 
performed in two steps; for the interpretation of these steps with respect to applications in radiometry
we refer to Helstrom's paper.

In the two steps described by Helstrom we need two given numbers $q_0, q_1$, satisfying  $0<q_0\le q_1<1$. For the asymptotic inversion we assume that 
$\mu$ is a large parameter. The two steps are:

\begin{description}
\item[Step~1:]
Find $y$ from the  equation 
\begin{equation}\label{eq:inv01}
Q_\mu(0,y)=q_0,
\end{equation}
and denote this value by $y_0$. Recall, see \eqref{eq:def07},  that $Q_\mu(0,y)=Q_\mu(y)$ (the normalized incomplete gamma function).
\item[Step~2:]
Find $x$ from the  equation 
\begin{equation}\label{eq:inv02}
Q_\mu(x,y_0)=q_1,
\end{equation}
and denote this value by $x_1$. The value $y_0$ is obtained in Step~1.
\end{description}

For Step~1 we refer to \cite{Temme:1992:AIG}; see also \cite[Chapter~10]{Gil:2007:NSF} and \cite{Gil:2010:AIC}. For Step~2 we use an asymptotic representation of $Q_\mu(x,y)$, and for these details we refer to \S\ref{sec:asinvx}. This inversion process is with respect to $x$ with $y=y_0$ a fixed value obtained in Step~1. 

We also consider the inversion with respect to $y$ with fixed $x$; see \S\ref{sec:asinvy}.

In statistics, the inversion of  $Q_\mu(x,y)$ with respect to $x$ corresponds to 
the problem of inverting the distribution function with
respect to the non-centrality parameter given the upper tail probability.
On the other hand, the inversion of $P_{\mu}(x,y)$ with respect to $y$ with fixed $x$
corresponds to the problem of computing the $p$-quantiles of the distribution function. 

Before this, we discuss the convergence of iterative methods (Newton, secant). The asymptotic methods that
will be described in section \S\ref{sec:asinv} will  then be used to boost the convergence of iterative methods, or,  when $\mu$ is large enough, they can be used as the sole method of inversion.

\section{Inversion by iterative methods}
\label{iterative}

We discuss in detail the inversion with respect to $x$, that is, for Step~2. For the inversion with respect to $y$ similar results can be obtained.

First we observe that values $q_1\le q_0$ in Step~2 do not make sense. For an explanation, consider the quadrant $(x\ge0,y\ge0)$.
When we have found a value $y_0$ from Step~1 in the inversion process, satisfying $q_0=Q_\mu(0,y_0)=Q_\mu(y_0)$, the inversion of $Q_\mu(x,y_0)=q_1$ can be interpreted as the inversion with respect to $x$ in the $(x\ge0,y\ge0)$ quadrant along the horizontal line $y=y_0$. When we follow $Q_\mu(x,y_0)$ along this line, starting at $x=0$ with value $q_0$, observing that this function is increasing with increasing $x$, it can never be equal to $q_1$ when $q_1<q_0$. When $q_1=q_0$ the solution of the inversion is $x=0$.

Since, with $y_0$, fixed the function $Q_\mu(x,y_0)$ is increasing as a function of $x$, it changes from concave to convex.
This is a favorable situation to ensure convergence
of the Newton method in the second step of inversion, that is, for the computation of the value of $x$ (called $x_1$)
such that $Q_{\mu}(x,y_0)=q_1$, with $q_1>Q_{\mu}(y_0)$ .

Because the function $f(x)=Q_{\mu}(x,y_0)-q_1$ is such that $f'(x)>0$, $f''(x)<0$ if $x>x^{+}=y-\mu-\frac12$, simple graphical
arguments show that if $x_1>x^{+}$ the Newton iteration
\begin{equation}
\bar{x}_{n+1}=\bar{x}_n -\Frac{f(\bar{x}_n)}{f'(\bar{x}_n)}
\end{equation}
converges monotonically to $x_1$ for any starting value $\bar{x}_0\in [x^{+},x_1]$. Similarly, if $x_1<x^{-}=x-\mu-1$, 
monotonic convergence is guaranteed for any starting value $\bar{x}_0\in [x_1,x^{-}]$. 

Therefore, with $\bar{x}_0$ equal to $x^{+}$ and $x^{-}$ we have guaranteed convergence in the cases $x_1>x^{+}$ and
$x_1 <x^{-}$, respectively. 
On the other hand, in the interval $(x^-,x^+)$ the second derivative is small and the graph is
close to a straight line, in which case the Newton method is safe. For this reason, any starting value in $[x^-,x^+]$ 
produces convergence.  If $y-\mu-1<0$, we can take as starting value $\bar{x}_0=0$.

For the Newton method, it is necessary to compute the derivative of $f(x)$, which can be written in terms of the Bessel function
(\ref{eq:def10}). An alternative is the secant method, which has a slightly smaller convergence rate but avoids the use
of the derivative. Convergence also holds for the secant method, for similar graphical reasons as for Newton method,
and the choice of initial values $x^{-}$ and $x^{+}$ is a natural and safe selection; if these values are negative, we
take small positive values of $x$ as starting values.

Depending on the value of $q_1$, it may be more interesting to invert using the $P$-function instead of the $Q-$function.
If $q_1$ is close to $1$ it is more interesting to consider the value $p_1=1-q_1$ and invert the equation 
$P_{\mu}(x,y_0)=p_1$. Inversion for $P$ can also be carried out reliably with the Newton or the secant method 
and the same starting values used for $Q$ can be used for $P$ (monotonicity an convexity properties for $P$ immediately
follow from those of $Q$ considering (\ref{eq:intro03})).

Inversion with respect to $y$ can be done in a similar way. $Q(x,y)$ is decreasing as a function of $y$ and with
the convexity properties of Eq.~(\ref{conve_y}) we conclude that choosing $y^{-}=x+\mu-\frac32$ and $y^{+}=x+\mu-1$ is
a safe selection ensuring convergence (and if these values are negative,
taking instead small starting values for $y$).

Inversion using the secant method has been exhaustively tested and convergence holds in parameter
ranges where our algorithm for the Marcum-$Q$ function works \cite{Gil:2013:CMQ}; more that $10^8$ test points have
been considered.

\section{Asymptotic inversion}\label{sec:asinv}
For the asymptotic inversion methods we use the representation 
\begin{equation}\label{eq:asinv01}
Q_\mu(\mu x,\mu y)=\tfrac12\erfc\left(\zeta\sqrt{\mu/2}\right)-R_\mu(\zeta).
\end{equation}
Here, $\erfc\,z$ is the complementary error function defined by
\begin{equation}\label{eq:asinv02}
\erfc\,z=\frac{2}{\sqrt{\pi}}\int_z^\infty e^{-t^2}\,dt,
\end{equation}
and $R_\mu(\zeta)$ can be written in the form
\begin{equation}\label{eq:asinv03}
R_\mu(\zeta)=\frac{e^{-\frac12\mu\zeta^2}}{\sqrt{2\pi\mu}}S_\mu(\zeta), \quad 
S_\mu(\zeta)\sim \sum_{n=0}^\infty \frac{d_n(\zeta)}{\mu^n}.
\end{equation}
The first coefficient is
\begin{equation}\label{eq:asinv04}
d_0(\zeta)=\frac{1}{\zeta}-\frac{1+2x+\sqrt{1+4xy}}{2(y-x-1)(1+4xy)^{1/4}}.
\end{equation}
The quantity $\zeta$  is given by
\begin{equation}\label{eq:asinv05}
\tfrac12\zeta^2= x + y-\sqrt{1+4xy}+\ln\frac{1+\sqrt{1+4xy}}{2y},
\end{equation}
with $\sign(\zeta)=\sign(y-x-1)$.

For details on the representation in \eqref{eq:asinv01} we refer to Appendix~B. For the inversion process we consider the two cases: inversion with respect to $x$ and to~$y$.

\begin{remark}\label{rem:marc01}
{\rm
For the complementary function we have
\begin{equation}\label{eq:asinv06}
P_\mu(\mu x,\mu y)=\tfrac12 \erfc\left(-\zeta\sqrt{\mu/2}\right)
+R_\mu(\zeta),
\end{equation}
where $R_\mu(\zeta)$ and $S_\mu(\zeta)$ are the same, with expansion as in \eqref{eq:asinv03}.
}
\end{remark}

\begin{remark}\label{rem:marc02}
{\rm
For the incomplete gamma functions we have similar representations (for convenience some of the notation is as in \eqref{eq:asinv01})
\begin{equation}\label{eq:invone01}
\begin{array}{ll}
Q_\mu(\mu y)=\tfrac12\erfc\left(\eta\sqrt{\mu/2}\right)-{R}_\mu(\eta),\\[8pt]
P_\mu(\mu y)=\tfrac12\erfc\left(-\eta\sqrt{\mu/2}\right)+{R}_\mu(\eta),
\end{array}
\end{equation}
where 
\begin{equation}\label{eq:invone02}
\tfrac12\eta^2=y-1-\ln y,\quad \sign(\eta)=\sign(y-1).
\end{equation}
The function ${R}_\mu(\eta)$ has the asymptotic representation
\begin{equation}\label{eq:invone03}
{R}_\mu(\eta)=\frac{e^{-\frac12\mu\eta^2}}{\sqrt{2\pi\mu}}{S}_\mu(\eta), \quad {S}_\mu(\eta)\sim \sum_{n=0}^\infty \frac{C_n(\eta)}{\mu^n}, \quad C_0(\eta)=\frac{1}{y-1}-\frac{1}{\eta}.
\end{equation}
The expansion is valid uniformly with respect to  $y\ge0$. The coefficient $C_0(\eta)$ is regular at the transition point $y=1$, and the same for all higher coefficients. For more details we refer to \cite[Chapter~11]{Temme:1996:SFA}.
}
\end{remark}

\subsection{Asymptotic inversion with respect to \protectbold{x}}\label{sec:asinvx}
This corresponds to Step~2 described in \S\ref{sec:MQinv}.  Throughout this case we take $y=y_0$, the value obtained in Step~1, that satisfies $q_0=Q_\mu(0,y_0)=Q_\mu(y_0)$, and the present inversion problem is 
\begin{equation}\label{eq:asinvx01}
Q_\mu(x,y)=q_1, \quad q_1>q_0.
\end{equation}

We use the method described for the incomplete gamma functions; see \cite{Temme:1992:AIG}, \cite[Chapter 10]{Gil:2007:NSF}, \cite{Gil:2010:AIC} and \cite{Gil:2012:IGR}.
We use representation \eqref{eq:asinv01} and start with solving the equation
\begin{equation}\label{eq:asinvx02}
\tfrac12\erfc\left(\zeta\sqrt{\mu/2}\right)=q_1.
\end{equation}
We call the solution $\zeta_0$, and considering $q_1$ as a function of $\zeta_0$, we have
\begin{equation}\label{eq:asinvx03}
\frac{dq_1}{d\zeta_0}=-\sqrt{\frac{\mu}{2\pi}}\,e^{-\frac12\mu\zeta_0^2}.
\end{equation}
Considering $q_1$ in \eqref{eq:asinvx01} as a function of $\zeta$, we have (see \eqref{eq:def10})
\begin{equation}\label{eq:asinvx04}
\frac{dq_1}{d\zeta}=\Frac{\partial Q_\mu(x,y)}{\partial x}\Frac{\partial x}{\partial \zeta}=
\left(\Frac{y}{x}\right)^{\mu/2} e^{-x-y} I_{\mu} \left(2\sqrt{xy}\right)\Frac{\partial x}{\partial \zeta}.
\end{equation}
Upon dividing, and replacing $x,y$ with $\mu x, \mu y$, we obtain
\begin{equation}\label{eq:asinvx05}
\frac{d\zeta}{d\zeta_0}=-\Frac{1}{\sqrt{2\mu\pi}}\,\frac{e^{-\mu(\frac12\zeta_0^2-x-y+\ln\rho)}}{I_\mu(\mu\xi)}
\Frac{\partial \zeta}{\partial x},
\end{equation}
where
\begin{equation}\label{eq:asinvx06}
\rho=\sqrt{\frac{y}{x}},\quad \xi=2\sqrt{xy}.
\end{equation}

Next we use an asymptotic representation of  the modified Bessel function that is valid for large values of $\mu$, uniformly with respect to $\xi\ge0$. We have\footnote{http://dlmf.nist.gov/10.41.E3}
\begin{equation}\label{eq:asinvx07}
I_{{\mu}}(\mu \xi)=\frac{e^{{\mu\eta}}}{\sqrt{2\pi\mu}\,(1+\xi^{2})^{{\frac{1}{4}}}}
T_\mu(\xi), \quad T_\mu(\xi)\sim\sum _{{k=0}}^{\infty}\frac{U_{k}(p)}{\mu^{k}},
\end{equation}
where
\begin{equation}\label{eq:asinvx08}
\eta=\sqrt{1+\xi^2}+\ln\frac{\xi}{1+\sqrt{1+\xi^2}}, \quad p=\frac{1}{\sqrt{1+\xi^2}},
\end{equation}
and the coefficients $U_k(p)$ are polynomials in $p$. The first two are
\begin{equation}\label{eq:asinvx09}
U_0(p)=1,\quad U_1(p)=\tfrac{1}{24}\left(3p-5p^3\right).
\end{equation}

Using the first part of \eqref{eq:asinvx07} in \eqref{eq:asinvx05} we obtain
\begin{equation}\label{eq:asinvx10}
f(\zeta)\frac{d\zeta}{d\zeta_0}=e^{-\frac12\mu(\zeta_0^2-\zeta^2)},
\end{equation}
where
\begin{equation}\label{eq:asinvx11}
\tfrac12\zeta^2=x+y-\ln\rho-\eta,
\end{equation}
which is the same as the relation in \eqref{eq:asinv05}, and
\begin{equation}\label{eq:asinvx12}
f(\zeta)=-\frac{T_\mu(\xi)}{\Frac{\partial \zeta}{\partial x}(1+\xi^{2})^{{\frac{1}{4}}}},
\end{equation}
in which 
\begin{equation}\label{eq:asinvx13}
\Frac{\partial \zeta}{\partial x}= \frac{1+\sqrt{1+\xi^2} -2y}{\zeta\left(1+\sqrt{1+\xi^2}\right)}.
\end{equation}

As in the asymptotic inversion of the incomplete gamma functions, where we obtained a similar equation as in \eqref{eq:asinvx10}, we solve this differential equation by substituting an expansion of the form
\begin{equation}\label{eq:asinvx14}
\zeta\sim\zeta_0+\sum_{n=1}^\infty\frac{\zeta_n}{\mu^n},
\end{equation}
where we have $\zeta_0$ computed from equation \eqref{eq:asinvx02}. After substituting this expansion in \eqref{eq:asinvx10} and considering the first order of approximation for large $\mu$, we obtain for $\zeta_1$ the relation
\begin{equation}\label{eq:asinvx15}
f(\zeta_0)=e^{\zeta_0\zeta_1}  \quad \Longrightarrow \quad \zeta_1=\frac{1}{\zeta_0}\ln\left(f(\zeta_0)\right),
\end{equation}
where
\begin{equation}\label{eq:asinvx16}
f(\zeta_0)=
-\frac{\zeta_0\left(1+\sqrt{1+\xi^2}\right)}{\left(1+\sqrt{1+\xi^2} -2y\right)(1+\xi^{2})^{{\frac{1}{4}}}},
\end{equation}
which can be written in the form
\begin{equation}\label{eq:asinvx17}
f(\zeta_0)=
\frac{\zeta_0}{y-x-1}
\frac{1+2x+\sqrt{1+4xy}}{2(1+4xy)^{\frac{1}{4}}}.
\end{equation}
Observe that this quantity is related to the coefficient $d_0(\zeta)$ defined in \eqref{eq:asinv04}. We have
$f(\zeta_0)=1-\zeta_0d_0(\zeta_0)$.

In  equation \eqref{eq:asinvx17} we may write $x$ also with subscript  $0$, because it is related to $\zeta_0$ via equation \eqref{eq:asinv05} with $\zeta$ replaced by $\zeta_0$. We have to invert this equation to find $x$ when $\zeta_0$ is available. This can be done by using standard equation solvers, such as Newton's method; see \S\ref{sec:moreinv} for more details.

When we have the value $x$  corresponding to $\zeta_0$, we can compute $f(\zeta_0)$ by using 
\eqref{eq:asinvx17}, and then $\zeta_1$ from \eqref{eq:asinvx15}.  This gives the the second-order approximation  $\zeta\sim \zeta_0+\zeta_1/\mu$. When we have higher-order coefficients $\zeta_n$ in the expansion of $\zeta$, we use this $\zeta$ to find $x$ from \eqref{eq:asinv05}, and so on. These higher coefficients can be obtained by
expanding $f(\zeta)$ in negative powers of $\mu$ (after substitution of \eqref{eq:asinvx14}, and also by expanding the exponential function in \eqref{eq:asinvx10}). The comparison of the coefficients of equal powers of $\mu$ gives the relations for the $\zeta_j$.

\subsection{Asymptotic inversion with respect to \protectbold{y}}\label{sec:asinvy}

In this case the inversion problem is defined by

\begin{equation}\label{eq:asinvy01}
Q_\mu(x,y)=q,\quad 0<q<1,
\end{equation}
with $x$ a given fixed value. Observe that this time $q$ may be any value in the interval $(0,1)$, because for any positive $x$, we have have $Q_\mu(x,0)=1$ and $Q_\mu(x,y)$ is monotonically decreasing to $0$ as $y\to\infty$; see \eqref{eq:def07} and \eqref{eq:def10}.

We proceed as in the previous case, computing $\zeta_0$ from  equation \eqref{eq:asinvx02} with $q=q_1$, and by using \eqref{eq:def10}. We obtain the relation (cf.~\eqref{eq:asinvx05})
\begin{equation}\label{eq:asinvy02}
\frac{d\zeta}{d\zeta_0}=\rho\Frac{1}{\sqrt{2\mu\pi}}\,\frac{e^{-\mu(\frac12\zeta_0^2-x-y+\ln\rho)}}{I_{\mu-1}(\mu\xi)}
\Frac{\partial \zeta}{\partial y}.
\end{equation}

We replace the Bessel function by using $I_{\mu-1}(\mu\xi)=I_\mu^\prime(\mu\xi)+(1/\xi)I_\mu(\mu\xi)$ and we use the expansion for the derivative
\begin{equation}\label{eq:asinvy03}
I_{{\mu}}^\prime(\mu \xi)=\frac{e^{\mu\eta}}{\xi \sqrt{2\pi\mu}}\,(1+\xi^{2})^{\frac{1}{4}}
W_\mu(\xi), \quad W_\mu(\xi)\sim\sum _{{k=0}}^{\infty}\frac{V_{k}(p)}{\mu^{k}},
\end{equation}
where $\eta$ and $p$  are the same as in \eqref{eq:asinvx08} and  the coefficients $V_k(p)$ are polynomials in $p$. The first two are
\begin{equation}\label{eq:asinvy04}
V_0(p)=1,\quad V_1(p)=\tfrac{1}{24}\left(-9p+7p^3\right).
\end{equation}

This gives the analogue of \eqref{eq:asinvx10} in the form
\begin{equation}\label{eq:asinvy05}
g(\zeta)\frac{d\zeta}{d\zeta_0}=e^{-\frac12\mu(\zeta_0^2-\zeta^2)},
\end{equation}
where
\begin{equation}\label{eq:asinvy06}
g(\zeta)=\frac{T_\mu(\xi)+\sqrt{1+\xi^2}W_\mu(\xi)}{2y\Frac{\partial \zeta}{\partial y}(1+\xi^{2})^{{\frac{1}{4}}}},
\end{equation}
in which 
\begin{equation}\label{eq:asinvy07}
\Frac{\partial \zeta}{\partial y}= \frac{y-2xy-1+(y-1)\sqrt{1+\xi^2}}{y\zeta\left(1+\sqrt{1+\xi^2}\right)}.
\end{equation}

When we take $T_\mu(\xi)=1$ and $W_\mu(\xi)=1$ we obtain $g_0(\zeta_0)$, which turns out to be the same as $f(\zeta_0)$ given in \eqref{eq:asinvx17}.
In this way we obtain for $\zeta$ the approximation $\zeta\sim\zeta_0+\zeta_1/\mu$, where $\zeta_1$ is given in \eqref{eq:asinvx15}.

\begin{example}\label{ex:trans}
{\rm
As an application we use $q=\frac12$. Then the equation
\begin{equation}\label{eq:asinvy12}
\tfrac12\erfc\left(\zeta_0\sqrt{\mu/2}\right)=\tfrac12
\end{equation}
gives $\zeta_0=0$ and we have
\begin{equation}\label{eq:asinvy13}
\zeta\sim \frac{\zeta_1}{\mu}= \frac{d_0(x)}{\mu},
\end{equation}
where $d_0(x)$ is shown in the expansion \eqref{eq:moreinv09},
and using this in expansion  \eqref{eq:moreinv05} we find
\begin{equation}\label{eq:asinvy14}
y\sim x+1+\frac{b_1(x)d_0(x)}{\mu}=x+1-\frac{3x+1}{3\mu(2x+1)},\quad \mu\to \infty.
\end{equation}

This is for the scaled variables in $Q_\mu(\mu x,\mu y)$.  For the real life variables $x,y$ in $Q_\mu(x,y)$ we have: $Q_\mu(x,y)=\frac12$ when
\begin{equation}\label{eq:asinvy15}
y\sim x+\mu-\frac{3x+1}{3(2x+1)},\quad \mu\to \infty.
\end{equation}
This gives a description of the transition  in the quadrant $(x\ge0,y\ge0)$ from small values to values near unity of the Marcum functions.
}
\end{example}

\begin{remark}\label{rem:pinversion}
{\rm
For applications in mathematical statistics it is of interest to consider the inversion for the $P-$function in the form
\begin{equation}\label{eq:asinvy16}
P_\mu(x,y)=p,\quad 0<p<1.
\end{equation}
By using $q=1-p$ we can use the inversion of the $Q-$function, but when $p$ is very small the evaluation $q=1-p$ does not make sense. We can repeat the analysis for the $P-$function, and the only change we have to make is to change the sign of $\zeta_0$, and to assume as earlier $\sign(\zeta_0)=\sign(y-x-1)$. This follows from the representation in Remark~\ref{rem:marc02}, see \eqref{eq:asinv06}.
}
\end{remark}

In Appendix~A we give more details on the inversion process.

\section{Numerical examples}\label{sec:numex}

Again we consider two cases: inversion with respect to $x$ and to $y$. We describe the inversion using 
the asymptotic methods. As we will see, we obtain good accuracy even for relatively small values of $\mu$.
If more accuracy is needed for small $\mu$, it is always possible to do inversion using directly the secant
method for the Marcum function, which can be computed using the algorithm in \cite{Gil:2013:CMQ}, 
as done in \S\ref{iterative}; the performance of the secant method is improved by using starting values provided
by the asymptotic methods when these are accurate.

For the inversion with respect to $x$ described in \S\ref{sec:asinvx}, we compare our method with 
the numerical results shown in \cite{Helstrom:1998:AIM}. In that paper the secant method 
has been used for finding the value of $x$ when $y$ and $\mu$ are given. Rather large values of $\mu$
and rather small values 
of $q_0$ and $1-q_1$ are considered. From the point of view of the application in radiometry considered
in \cite{Helstrom:1998:AIM}, this
is due to the fact that $q_0$ represents a false-alarm probability (with typical low values) and $q_1$
represents a probability of detection (with typical large values).  
We use the same values of that paper $q_0=\left\{10^{-6},10^{-8}\right\},\,q_1=\left\{0.9,\,0.999\right\}$, and one other set of $q_0$ and $q_1$.

In Table~\ref{tab:numer01} we give the relative errors of the computations.
For  Step~1 of the inversion process the value $y_0$ satisfying $Q_\mu(y_0)=q_0$ is computed, 
and with this $y_0$ the value $\widetilde{q}_0=Q_\mu(y_0)$ by using the algorithms for 
the incomplete gamma function ratios described in \cite{Gil:2012:IGR}. The  displayed relative
 error is  $\delta_0=\vert q_0/\widetilde{q}_0-1\vert$.

In Step~2 we have computed  the value $x_1$ satisfying $Q_\mu(x_1,y_0)=q_1$, and computed $\widetilde{q}_1=Q_\mu(x_1,y_0)$ by
 using the algorithms for the Marcum functions described in \cite{Gil:2013:CMQ}.
Then, the displayed relative error is  $\delta_1=\vert q_1/\widetilde{q}_1-1\vert$.
To test  the inversion values for $\mu=10,\,20,\,50,\,100,\,200,\,500,\,1000$,
we have used the fixed-precision Fortran 90 module {\bf MarcumQ} presented in \cite{Gil:2013:CMQ}. The admissible range of computation in {\bf MarcumQ} is restricted to $\mu \le 10^4$, otherwise  overflow/underflow problems in IEEE double-precision arithmetic may occur.
An extended precision version of the algorithms has been used to test the values $\mu=10^5,\,10^7,\,10^9$.
 
For the inversion in Step~2 we have used only the approximations  
obtained from the asymptotic inversion process described in \S\ref{sec:asinvx}, and not the algorithms for the Marcum function  $Q_\mu(x,y)$ themselves.
The relative errors in Table~\ref{tab:numer01} are much smaller than those given in \cite{Helstrom:1998:AIM}, which proves the
accuracy of our asymptotic inversion method in Step~2. The results for Step~1 are also better, but for that case we have used not only asymptotic methods, but highly-accurate numerical algorithms from \cite{Gil:2012:IGR}.  

\begin{figure}
\begin{center}
\epsfxsize=8cm \epsfbox{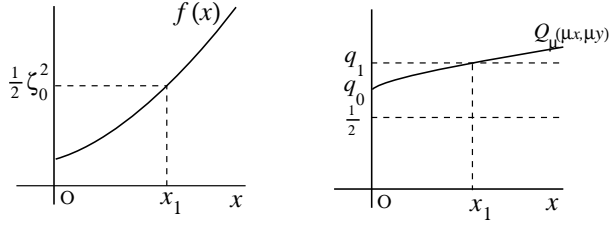}
\caption{The function $f(x)$ defined in \eqref{eq:moreinvx01}, $x\ge0$, $y<1$ (left) and the inversion with  $q_1>\frac12$ (right).  \label{fig:plot1}}
\end{center}
\end{figure}

For applying Newton's method in Step~2 we have to use different starting choices.
\begin{enumerate}
\item
When $y<1$ (see \S\ref{sec:moreinvx1} and Figure~\ref{fig:plot1}) we fit a polynomial with the values $f(0)$, $f^\prime(0)$ and $f^{\prime\prime}(0)$.

\item
When $y\ge1$ and $q_1<\frac12$ (see \S\ref{sec:moreinvx2} and Figure~\ref{fig:plot2} (left and middle)) we need to compute the $x_1$ left of $y-1$,  and start at $x=0$, without further~fitting.
\item
When $y\ge1$ and $q_1>\frac12$ (see  Figure~\ref{fig:plot2} (left and right)) we need the $x_1$ on the right of $y-1$. Observe that $f(x)>g(x)= x+y-\sqrt{1+4xy}$ if $x>y-1$. The function $g(x)$ has zeros at $x=y-1$  
and $x=y+1$ (and between these values it is negative). A suitable starting value is the solution of the equation $g(x)=\frac12\zeta^2$, which root is located on the right of $x=y+1$.
\end{enumerate}

With the values of $\zeta_0$,  $x=x_1$ and $y=y_0$ we can compute $f_0(\zeta_0)$ given in \eqref{eq:asinvx17}, and then $\zeta_1$ of \eqref{eq:asinvx15}, either using their explicit expressions or the expansions in \eqref{eq:moreinv07}.
For the numerical examples shown in Table~\ref{tab:numer01} we have used these series expansions with terms up to $k=10$ when needed.

Then we can compute $\zeta\sim\zeta_0+\zeta_1/\mu$ (see \eqref{eq:asinvx14}), and using this value of $\zeta$ in \eqref{eq:moreinv01} yields a new value $x$. In the described inversion process and in the computations, the $x$ and $y$ are scaled variables. The real life $x$ and $y$  are replaced with  $\mu x$ and $\mu y$, and the scaled values solve the equation  $Q_\mu(\mu x,\mu y)=q_1$.

\renewcommand{\arraystretch}{1.2}
\begin{table}
  \caption{Relative errors for the inversions in Step~1 and Step~2. For the meaning of the relative errors $\delta_0$ and $\delta_1$ we refer to the text of \S\ref{sec:numex}.
  \smallskip
  \label{tab:numer01}}
  \begin{centering}
  {
  {
    \begin{tabular}{c|cc|cc|cc}\hline\smallskip
 & $q_0$=1.0e-6 & $q_1$=0.9  & $q_0$=1.0e-8 & $q_1$=0.999 & $q_0$=0.4 & $q_1$=0.6 \\
\cline{2-7}\smallskip 
$\mu$ & $\delta_0$ & $\delta_1$  & $\delta_0$ & $\delta_1$ & $\delta_0$ & $\delta_1$ \\
\hline
 1.00e+1   & 1.27e-15  & 1.23e-05   &  2.31e-15 & 2.30e-07   & 1.39e-16  & 1.16e-04   \\
 2.00e+1   & 1.48e-15  & 7.43e-06   &  0        & 7.29e-08   & 1.39e-16  & 1.94e-05   \\
 5.00e+1   & 1.48e-15  & 3.27e-06   &  1.48e-15    & 3.70e-08   & 1.39e-16  & 3.09e-06   \\
 1.00e+2   & 8.68e-15  & 1.55e-06   &  8.93e-15    & 1.97e-08   & 1.39e-16  & 3.15e-06   \\
 2.00e+2   & 2.54e-15  & 6.51e-07   &  2.98e-15    & 9.34e-09   & 5.55e-16  & 1.76e-06   \\
 5.00e+2   & 9.11e-15  & 1.68e-07   &  1.67e-14    & 2.90e-09   & 1.39e-16  & 6.18e-07   \\
 1.00e+3   & 1.55e-14  & 4.98e-08   &  1.85e-14    & 1.02e-09   & 0  & 2.54e-07   \\
 1.00e+5   & 2.10e-17  & 7.14e-11  & 3.81e-17  & 1.01e-12  & 1.02e-17  & 3.41e-10  \\
 1.00e+7   & 2.96e-16  & 1.04e-13  & 6.25e-16  & 1.75e-15  & 1.37e-16  & 3.51e-13  \\
 1.00e+9   & 1.30e-15  & 1.66e-13  & 7.16e-18  & 7.94e-17  & 3.16e-16  & 7.90e-16  \\
\hline
    \end{tabular}} \\}
  \end{centering}
  \end{table}
\renewcommand{\arraystretch}{1.0}

For the inversion with respect to $y$ we apply the method described in \S\ref{sec:asinvy}, which  is straightforward. We show in Table~\ref{tab:numer02}
the results of numerical computations for several values of $\mu$ and in each case we take $x=\mu$. We take in the approximation $\zeta\sim\zeta_0+\zeta_1/\mu$ in the columns indicated by $\zeta_0$, $\delta_0$ only the term $\zeta_0$,  and in the columns indicated by $\zeta_1$, $\delta_1$, the approximation including $\zeta_1$.  The latter case has a better performance, as expected.

The values $\delta_0$, $\delta_1$ are the corresponding relative errors $\vert q/\widetilde q -1\vert$, where $q$ is the given value and $\widetilde q=Q_\mu(x,y)$, with the given $\mu$ and $x$, and the $y$ obtained by inversion.

\renewcommand{\arraystretch}{1.2}
\begin{table}
  \caption{Relative errors for the inversion values in the inversion process with respect to $y$. For the meaning of the 
$\zeta_0$, $\delta_0$ and $\zeta_1$, $\delta_1$ we refer to the text of \S\ref{sec:numex}.
  \smallskip
  \label{tab:numer02}}
  \begin{centering}
  {
  {
    \begin{tabular}{c|cc|cc|cc}\hline\smallskip
  & \ $q$=1.0e-6 &   & \ $q$=0.5 & &\  $q$=0.9999 &  \\
\cline{2-7}\smallskip 
$\mu$ & $\zeta_0$,\ $\delta_0$ & $\zeta_1$,\ $\delta_1$  & $\zeta_0$,\ $\delta_0$ & $\zeta_1$,\ $\delta_1$ & $\zeta_0$,\ $\delta_0$ & $\zeta_1$,\ $\delta_1$ \\
\hline
1.00e+1   & 3.80e-01   & 1.89e-03   & 6.95e-02   & 5.65e-04   & 4.83e-05   & 3.92e-07   \\
 2.00e+1   & 2.72e-01   & 7.90e-04   & 4.81e-02   & 1.99e-04   & 3.00e-05   & 1.33e-07   \\
 5.00e+1   & 1.75e-01   & 2.34e-04   & 2.98e-02   & 5.01e-05   & 1.70e-05   & 3.07e-08   \\
 1.00e+2   & 1.24e-01   & 8.97e-05   & 2.09e-02   & 1.77e-05   & 1.14e-05   & 1.02e-08   \\
 2.00e+2   & 8.87e-02   & 3.35e-05   & 1.47e-02   & 6.25e-06   & 7.80e-06   & 3.45e-09   \\
 5.00e+2   & 5.63e-02   & 8.96e-06   & 9.24e-03   & 1.58e-06   & 4.78e-06   & 8.40e-10  \\
 1.00e+3   & 3.99e-02   & 3.25e-06   & 6.52e-03   & 5.59e-07   & 3.33e-06   & 2.91e-10  \\
  1.00e+5   & 4.00e-03   & 3.44e-09   & 6.47e-04   & 5.59e-10  & 3.22e-07   & 2.79e-13  \\
 1.00e+7   & 4.01e-04   & 3.46e-12  & 6.47e-05   & 5.59e-13  & 3.21e-08   & 2.74e-16  \\
 1.00e+9   & 4.02e-05   & 7.69e-14  & 6.47e-06   & 1.69e-15  & 3.21e-09   & 9.33e-17  \\
\hline
    \end{tabular}} \\}
  \end{centering}
  \end{table}
\renewcommand{\arraystretch}{1.0}

\bigskip

\section{Appendix~A: More details on the inversion process}\label{sec:moreinv}

We discuss details on the inversion of equation \eqref{eq:asinv05}
\begin{equation}\label{eq:moreinv01}
\tfrac12\zeta^2= x + y-\sqrt{1+4xy}+\ln\frac{1+\sqrt{1+4xy}}{2y},
\end{equation}
with $\sign(\zeta)=\sign(y-x-1)$ with respect to $x$ (with $y$ fixed) or to $y$ (with $x$ fixed). 

The value $\zeta$ follows from solving one of the equations (see \eqref{eq:asinvx02} and \eqref{eq:asinvy01}, respectively)
\begin{equation}\label{eq:moreinv02}
\tfrac12\erfc\left(\zeta\sqrt{\mu/2}\right)=q_1,\quad Q_\mu(x,y)=q,
\end{equation}
or from successive steps in the inversion process; see our remarks at the end of \S\ref{sec:asinvx}.

When $\zeta $ is small, $\vert y-x-1\vert$  is small, and it is convenient to have  expansions of $y-x-1$ in powers of $\zeta$. For the inversion with respect to $x$ we have 
\begin{equation}\label{eq:moreinv03}
x=y-1+\sum_{k=1}^\infty a_k(y) \zeta^k.
\end{equation}
From \eqref{eq:moreinv01} we find (taking into account the relation $\sign(\zeta)=\sign(y-x-1)$ and assuming $2y>1$)
\begin{equation}\label{eq:moreinv04}
a_1(y)=-\sqrt{2y-1},\quad a_2(y)=\frac{3y-1}{3(2y-1)},\quad a_3(y)=\frac{6y-1}{36(2y-1)^{5/2}}.
\end{equation}

For the inversion with respect to $y$ we use
\begin{equation}\label{eq:moreinv05}
y=x+1+\sum_{k=1}^\infty b_k(x) \zeta^k,
\end{equation}
with first coefficients
\begin{equation}\label{eq:moreinv06}
b_1(x)=\sqrt{2x+1},\quad b_2(x)=\frac{3x+1}{3(2x+1)},\quad b_3(x)=\frac{6x+1}{36(2x+1)^{5/2}}. 
\end{equation}

It is also convenient to have the following expansions
\begin{equation}\label{eq:moreinv07}
f(\zeta_0)=\sum_{k=0}^\infty c_k(x) \zeta_0^k, \quad
\zeta_1=\sum_{k=0}^\infty d_k(x) \zeta_0^k,
\end{equation}
where $f(\zeta_0)$ and $\zeta_1$ are defined in \eqref{eq:asinvx15} and \eqref{eq:asinvx16}.

The first coefficients are
\begin{equation}\label{eq:moreinv08}
\begin{array}{ll}
\dsp{c_0(x)=1, \quad c_1(x)=-\frac{3x+1}{3(2x+1)^{3/2}},\quad c_2(x)=\frac{18x^2+6x+1}{12(2x+1)^3},}\\[8pt]
\dsp{c_3(x)=-\frac{675x^3+81x^2+36x+4}{270(2x+1)^{9/2}},}
\end{array} 
\end{equation}
\begin{equation}\label{eq:moreinv09}
\begin{array}{ll}
\dsp{d_0(x)=-\frac{3x+1}{3(2x+1)^{3/2}},\quad d_1(x)=\frac{36x^2+6x+1}{36(2x+1)^3},}\\[8pt]
\dsp{d_2(x)=-\frac{2160x^3-594x^2-9x-1}{1620(2x+1)^{9/2}}.}
\end{array} 
\end{equation}

\subsection{Details on the inversion with respect to \protectbold{x}}\label{sec:moreinvx}

We give a few details of this inversion, because different cases have to be considered.
Let us denote the right-hand side of \eqref{eq:moreinv01} by $f(x)$, that is,
\begin{equation}\label{eq:moreinvx01}
f(x)=x + y-\sqrt{1+4xy}+\ln\frac{1+\sqrt{1+4xy}}{2y}.
 \end{equation}
Then 
\begin{equation}\label{eq:moreinvx02}
f^\prime(x)=\frac{1-2y+\sqrt{1+4xy}}{1+\sqrt{1+4xy}},
 \end{equation}
and both $f(x)$ and $f^\prime(x)$ vanish at $x =y-1$, with $f^{\prime\prime}(y-1)=1/(2y-1)$, this value being well defined because $x=y-1$ and $y=\frac12$ cannot happen, because then $x=-\frac12$ (a special case of the vanishing of $\sqrt{1+4xy}$\,).

\subsubsection{The case \protectbold{y<1}}\label{sec:moreinvx1}
For this case we refer to Figure~\ref{fig:plot1}. When $y<1$, $f(x)$ is monotonically increasing, starting with $f(0)=y-1-\ln(y)$, and the inversion of $\frac12\zeta_0^2=f(x)$ with respect to $x$ can be done straightforwardly. The only point to verify is whether indeed
$\frac12\zeta_0^2>f(0)$, otherwise there is no real positive root of the equation $\frac12\zeta_0^2=f(x)$.

To verify this point, observe that the present $\zeta_0$  follows from the first equation in \eqref{eq:moreinv02}, and because we assume that $y<1$, the values $q_0, q_1$ should satisfy $q_1>q_0>\frac12$; see Figure~\ref{fig:plot1}(right). This means that the inversion   
of  $Q_\mu(\mu x,\mu y)=q_1$ happens in the quadrant $(x\ge0,y\ge0)$ below the line $y=x+1$ (scaled variables) 
on the horizontal line $y=y_0$, on which no transition point can be found, on which $Q_\mu(\mu x,\mu y)$ is increasing, and on which $\zeta$ of representation \eqref{eq:asinv01} is increasing in absolute value (with $\zeta<0$).  The starting value of this $\zeta$ (at $x=0$) is the $\eta$ in the representation of the incomplete gamma function in \eqref{eq:invone01}, which satisfies $\frac12\eta^2=f(0)$. Hence, the $\zeta=\zeta_0$ in the relation \eqref{eq:moreinv01} corresponding to the $x$ and $y$ values satisfying $Q_\mu(\mu x,\mu y)=q_1$  is such that $\frac12\zeta_0^2>f(x)$, as shown in Figure~\ref{fig:plot1} (left). On the right we see the graph of $Q_\mu(\mu x,\mu y)$ with no transition point (inflection point).

\begin{figure}
\begin{center}
\epsfxsize=12cm \epsfbox{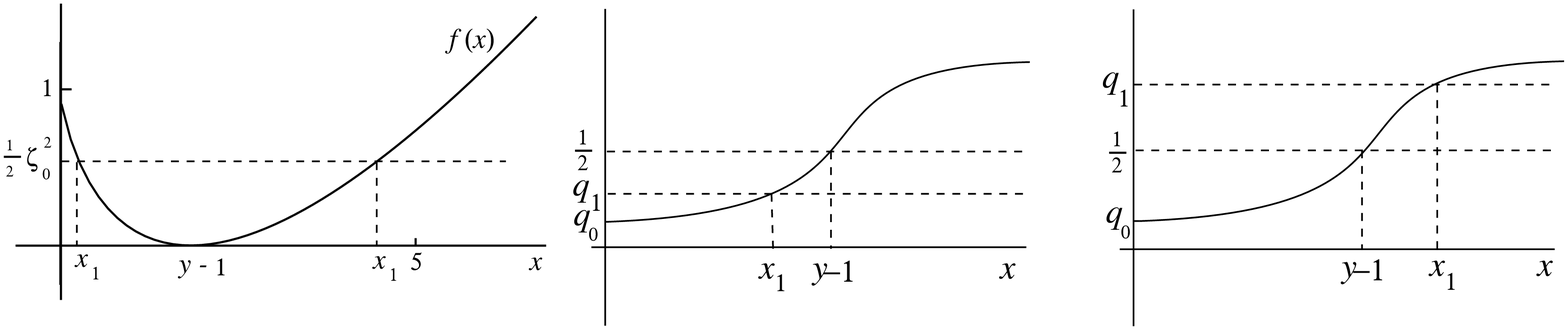}
\caption{The function $f(x)$ defined in \eqref{eq:moreinvx01}, $x\ge0$, $y=3$. The minimum occurs at $x=y-1=2$.  \label{fig:plot2}}
\end{center}
\end{figure} 

\subsubsection{The case \protectbold{y>1}}\label{sec:moreinvx2}
When $y>1$ the graph of $f(x)$ is as in Figure~\ref{fig:plot2} (left), where we see a minimum at $x=y-1$. There are two possible subs-cases. 
\begin{enumerate}
\item
$q_1<\frac12$. The function $Q_\mu(\mu x,\mu y)$ starts below $q_1$ at $(x,y)=(0,y_0)$, with positive value of $\zeta$ in representation \eqref{eq:asinv01} because $Q_\mu(\mu x,\mu y)< \frac12$ before $x$ crosses the value $y-1$, see Figure~\ref{fig:plot2} (middle). When $x$ increases, $\zeta$ becomes smaller until it becomes $0$ when $x=y-1$. As above, at $x=0$, $\zeta$ equals the corresponding value $\eta$ for the incomplete gamma function  in \eqref{eq:invone01}, with $\frac12\eta^2=f(0)$. This verifies, as shown in Figure~\ref{fig:plot2} (left), that $f(0)>\frac12\zeta_0$, and that the requested value $x_1$ is the one on the left in that figure.

 \item
$q_1>\frac12$. The function $Q_\mu(\mu x,\mu y)$ starts in Figure~\ref{fig:plot2} (right) at $q_0<\frac12$  (because $y>1$). As $x$ increases, $Q_\mu(\mu x,\mu y)$ becomes $\frac12$ (near $x=y-1$), then $Q_\mu(\mu x,\mu y)$ becomes larger than $\frac12$, below the line $y=x+1$. In this case the requested value $x_1$ satisfies $x_1>y-1$; it is the right one in Figure~\ref{fig:plot2} (left). The function $f(x)$ is monotonically increasing for $x>y-1$, and inversion is straightforward.

\end{enumerate}

\subsection{Details on the inversion with respect to \protectbold{y}}\label{sec:moreinvy}

In this case we solve the second equation in \eqref{eq:moreinv02} with respect to $y$, with given $x$. For small values of $\zeta$ and $\zeta_0$ we use the expansion given in \eqref{eq:moreinv02} 
and for other values a Newton process. The function on the right-hand side of \eqref{eq:moreinv01} (let us denote it by $g(y)$) has the derivative
\begin{equation}\label{eq:moreinvy01}
g^\prime(y)=\frac{y-2xy-1+(y-1)\sqrt{1+4xy}}{y\left(1+\sqrt{1+4xy}\right)}.
 \end{equation}
It vanishes with $g$ at $y=x+1$ and  $g$ is a convex function of $y$. For small values of $\zeta$ we can use the series in \eqref{eq:moreinv05} and for larger values Newton's method.

\section{Appendix~B: Asymptotic representation of\\
 the function \protectbold{Q_\mu(x,y)}}\label{sec:mqurep}

The asymptotic representation $Q_\mu(\mu x,\mu y)$ in \S\ref{sec:asinv} can be derived from the contour integral \begin{equation}\label{eq:mqurep03}
Q_\mu(\mu x,\mu y)=\frac{e^{-\mu\phi(\rho)}}{2\pi i} \int_{c-i\infty}^{c+i\infty}e^{\mu\phi(t)}\,\frac{dt}{\rho-t},\quad 0<c<\rho,
\end{equation}
where
\begin{equation}\label{eq:mqurep04}
\phi(t)=\tfrac12\xi\left(t+\frac{1}{t}\right)-\ln t,\quad \rho=\sqrt{\frac{y}{x}},\quad \xi=2\sqrt{xy}.
\end{equation}
For $P_\mu(\mu x,\mu y)$ a similar representation is valid when we take $c>\rho$. A slightly different form was derived in \cite[\S4]{Temme:1993:ANA}, and for details we refer to this paper.

In \eqref{eq:mqurep03} we can take $c=t_0$, where $t_0$ the positive saddle point of $\phi(t)$, which is given by $t_0=\left(1+\sqrt{1+\xi^2}\right)/\xi$. The saddle point coalesces with the pole at $\rho$ when $y=x+1$, and we assume that $0<t_0<\rho$. This implies  $y>x+1$, the domain in the $(x,y)-$plane where $Q_\mu(\mu x,\mu y)\le P_\mu(\mu x,\mu y)$ (approximately, but true for large $\mu$, $x$ and $y$).

The path of steepest descent $\calC$ through $t_0$ follows from the equation $\Im\phi(t)=0$. Let $t=r e^{i\theta}$, then we can describe $\calC$ by 
\begin{equation}\label{eq:mqurep06}
r=\frac{\theta}{\xi\sin\theta}+\sqrt{1+\frac{\theta^2}{\xi^2\sin^2\theta}},\quad -\pi<\theta<\pi.
\end{equation}

The transformation
\begin{equation}\label{eq:mqurep07}
\tfrac12 s^2 =\phi(t)-\phi(t_0)
\end{equation}
maps $\calC$ onto the imaginary axis in the $s-$plane. When taking the square root in this relation we assume $\sign(s)=\sign(t-t_0)$. The pole at $t=\rho$ corresponds to a pole in the $s-$plane at $s=\zeta$, say, where 
\begin{equation}\label{eq:mqurep08}
\tfrac12 \zeta^2 =\phi(\rho)-\phi(t_0).
\end{equation}
When taking the square root, we assume $\sign(\zeta)=\sign(\rho-t_0)=\sign(y-x-1)$ if $t>0$ and by continuity elsewhere. 

The transformation \eqref{eq:mqurep06} gives
\begin{equation}\label{eq:mqurep09}
Q_\mu(\mu x,\mu y)=\frac{e^{-\frac12\mu\zeta^2}}{2\pi i}\int_{-i\infty}^{i\infty}e^{\frac12\mu s^2}
f(s) \frac{ds}{\zeta-s},\quad f(s)= \frac{\zeta-s}{\rho -t}\frac{dt}{ds}.
\end{equation}
We split off the pole by writing $f(s)=\left(f(s)-f(\zeta)\right)+f(\zeta)$ and obtain
\begin{equation}\label{eq:mqurep10}
Q_\mu(\mu x,\mu y)=\tfrac12\erfc\left(\zeta\sqrt{\mu/2}\right)-\frac{e^{-\frac12\mu\zeta^2}}{2\pi i}\int_{-i\infty}^{i\infty}e^{\frac12\mu s^2}
g(s)\,ds,
\end{equation}
because $f(\zeta)=1$ and 
\begin{equation}\label{eq:mqurep12}
-\frac{e^{-\frac12\mu\zeta^2}}{2\pi i}\int_{-\infty}^{\infty}e^{-\frac12\mu \sigma^2}
 \frac{d\sigma}{\sigma+i\zeta}=\tfrac12\erfc\left(\zeta\sqrt{\mu/2}\right).
 \end{equation}
The function $g(s)$ is defined by $g(s)=(f(s)-f(\zeta))/(s-\zeta)$.
In \eqref{eq:mqurep10} we can drop the assumption $0<t_0<\rho$.

This produces the representation given in \S\ref{sec:asinvx}, and the coefficients $d_n(\zeta)$ in the expansion of $S_\mu(\zeta)$ follow from 
\begin{equation}\label{eq:mqurep13}
g(s)=\sum_{k=0}^\infty g_n s^n,\quad d_n(\zeta)= (-1)^n 2^n\left(\tfrac12\right)_n  g_{2n},
\end{equation}
where $\left(\frac12\right)_n=\Gamma\left(n+\frac12\right)/\Gamma\left(\frac12\right)$ (Pochhammer's symbol).

To compute the first coefficient 
\begin{equation}\label{eq:mqurep14}
d_0(\zeta)=g_0=g(0)=\frac{1-f(0)}{\zeta},\quad f(0)=\frac{\zeta}{\rho-t_0}\left.\frac{dt}{ds}\right\vert_{s=0},
\end{equation}
we first use
\begin{equation}\label{eq:mqurep15}
\frac{dt}{ds}=\frac{s}{\phi^\prime(t)}, \quad  \left.\frac{dt}{ds}\right\vert_{s=0}=\frac{t_0}{\left(1+\xi^2\right)^{1/4}},
\end{equation}
where we have used  l'H{\^o}pital's rule. The sign of the derivative at $s=0$ follows from the condition on the transformation in \eqref{eq:mqurep07}. This gives the coefficient $d_0(\zeta)$ shown in \eqref{eq:asinv04}.

\section*{Acknowledgements}
This work was supported by  {\emph{Ministerio de Econom\'{\i}a y Competitividad}}, 
project MTM2012-34787. The authors thank the referee for comments on the first version of the paper. NMT thanks CWI, Amsterdam, for scientific support.

\bibliographystyle{plain}

\bibliography{Marcum}

\end{document}